# Geometric Gamma Max-Infinitely Divisible Models

## Satheesh S

NEELOLPALAM, S. N. Park Road
Trichur – 680 004, **India.**
ssatheesh1963@yahoo.co.in

## Sandhya E

Department of Statistics, Prajyoti Niketan College
Pudukkad, Trichur – 680 301, **India.**
esandhya@hotmail.com

**Abstract.** A transformation of gamma max-infinitely divisible laws *viz.* geometric gamma max-infinitely divisible laws is considered in this paper. Some of its distributional and divisibility properties are discussed and a random time changed extremal process corresponding to this distribution is presented. A new kind of invariance (stability) under geometric maxima is proved and a max-AR(1) model corresponding to it is also discussed.

**AMS (2000) Classifications:** 60E05, 60F05, 60G70, 62M10.

**Keywords and phrases.** Max-infinite divisibility, geometric-max-infinite divisibility extremal processes, max-autoregressive process.

**1. Introduction.** Parallel to the classical notions of infinitely divisible (ID) laws and its subclass geometrically ID (GID) laws we have max-infinitely divisible (MID) and geometric-MID (G-MID) laws in the maximum setup which are discussed in Balkema and Resnick (1977), Rachev and Resnick (1991), Mohan (1998) and Satheesh (2002). For *d.f*s $F(x) = e^{-\psi(x)}$ that are MID (which is always true in **R**) distributions with *d.f*s of the form $\dfrac{1}{(1+\psi(x))}$ are referred to as G-MID laws. Processes related to these distributions are extremal processes, Rachev and Resnick (1991), Pancheva, *et al*. (200) and max-AR(1) processes, Satheesh and Sandhya (2006). Satheesh (2002) introduced φ-MID laws with *d.f* $\varphi\{-\log F(x)\}$ for a Laplace transform (LT) $\varphi$ and a *d.f F*. φ-MID laws can also be seen as the *d.f* obtained by randomizing the parameter $c>0$ in the first Lehman alternative $F^c$ obtained from *F*, by a distribution with LT $\varphi$. Setting $-\log F(x) = \psi(x)$ and taking $\varphi$ to be the LT of a gamma($\beta$) law we get the *d.f* $\dfrac{1}{(1+\psi(x))^\beta}$ which we will refer to as the *d.f* of a gamma-MID law. When *F* is max-semi-stable and $\varphi$ is exponential we get exponential max-semi-stable laws characterized in a max-AR(1) set up in Satheesh and Sandhya (2006). A gamma-max-semi-stable law was also discussed therein to illustrate the derivation of max-semi-selfdecomposable laws.

Dated 30 November 2007.

Pillai (1990) had introduced geometric exponential (G-exponential) laws having Laplace transform (LT) $\dfrac{1}{1+\log(1+\lambda)}$. Similarly, the LT $\dfrac{1}{1+\log(1+\lambda)^\beta}$, $\beta$>0 will be called as that of geometric gamma($\beta$) (G-gamma($\beta$)) law. Motivated by this construction, and writing the *d.f* of gamma-MID($\beta$) laws as $e^{-\beta\log(1+\psi(x))}$, $\beta$>0, we get the *d.f* of geometric-gamma-MID (G-gamma-MID($\beta$)) laws as $\dfrac{1}{\{1+\beta\log(1+\psi(x))\}}$. Geometric generalized gamma laws were introduced and studied in Sandhya and Satheesh (2007).

The purpose of this note is to discuss certain properties of gamma-MID and G-gamma-MID models. Potential applications of these models are in finance, insurance and stock market; see Kaufman (2001), Rachev (1993) and Mittnick and Rachev (1993). Distributional properties of these models including a new kind of invariance (stability) under geometric($p$)-maximum, are presented in section.2. In section.3 we discuss extremal processes and a max-AR(1) model related to them. The support of the distributions discussed here is **R**, and by a geometric($p$) law we mean a geometric law on $\{1,2,\ldots\}$ with mean $a = \frac{1}{p}$.

## 2. Divisibility properties of gamma-MID and G-gamma-MID laws.

With the above terminologies we have:

**Theorem.2.1** G-gamma-MID laws are G-MID.

*Proof.* We know that a *d.f* $F(x)$ is G-MID *iff* $e^{-\{\frac{1}{F(x)}-1\}}$ is MID, Rachev and Resnick (1991). From our construction we have the *d.f* $e^{-\beta\log(1+\psi(x))}$, $\beta$>0 that is always MID.

Setting $\dfrac{1}{F(x)}-1 = \beta\log\{1+\psi(x)\}$ we have $\dfrac{1}{1+\beta\log(1+\psi(x))} = F(x)$ is G-MID.

**Theorem.2.2** A *d.f* $F(x) = \dfrac{1}{1+\log(1+\psi(x))^\beta}$, $\beta$>0 is G-gamma-MID *iff* $\dfrac{1}{(1+\psi(x))^\beta}$, $\beta > 0$ is a *d.f.* This is clear from their construction above.

**Remark.2.1** It is interesting to note that $\dfrac{1}{(1-\log F(x))} = G(x)$ is a *d.f* for a given *d.f* $F(x)$ and in turn $\dfrac{1}{(1-\log G(x))}$ is also a *d.f* and this operation can be repeated to obtain new *d.f*s.

**Theorem.2.3** Every G-gamma-MID($\beta$) distribution is the limit distribution of geometric $\left(\frac{1}{n}\right)$-max of *i.i.d* gamma-MID $\left(\frac{\beta}{n}\right)$ variables as $n \to \infty$.

*Proof.* Let $F_n(x)$ denote the *d.f* of a geometric $\left(\frac{1}{n}\right)$-max of *i.i.d* gamma-MID $\left(\frac{\beta}{n}\right)$ variables.

Thus, $\quad F_n(x) = \dfrac{\frac{1}{n}(1+\psi(x))^{-\frac{\beta}{n}}}{1-\frac{n-1}{n}(1+\psi(x))^{-\frac{\beta}{n}}}$

$$= \dfrac{1}{n(1+\psi(x))^{\frac{\beta}{n}}-(n-1)}$$

$$= \dfrac{1}{1+n\{(1+\psi(x))^{\frac{\beta}{n}}-1\}}. \text{ Hence,}$$

$$\lim_{n\to\infty} F_n(x) = \dfrac{1}{1+\log(1+\psi(x))^{\beta}} = \dfrac{1}{1+\beta\log(1+\psi(x))}, \text{ proving the assertion.}$$

**Theorem.2.4** The limit of *n*-max of G-gamma-MID $\left(\frac{\beta}{n}\right)$ laws is gamma-MID($\beta$) as $n\to\infty$.

*Proof.* Since $\lim\limits_{n\to\infty}\left\{\dfrac{1}{1+\frac{\beta}{n}\log(1+\psi(x))}\right\}^n = e^{-\beta\log(1+\psi(x))} = \dfrac{1}{(1+\psi(x))^{\beta}}$, the claim is proved.

**Theorem.2.5** A distribution is invariant under geometric($p$)-max up to a scale-change *iff* it is geometric max-semi-stable with *d.f* $\dfrac{1}{1+\psi(x)} = \dfrac{1}{1+a\psi(bx)}$, $a>1$ and $b\in(0,1)\cup(1,\infty)$.

*Proof.* See theorem.3.2 in Satheesh and Sandhya (2006) and its proof which is formulated in the max-AR(1) set up and under the terminology of exponential-max-semi-stable law. If $b>1$ the geometric-max-semi-stable law is Frechet type and if $b<1$ it is Weibull type, see Satheesh and Sandhya (2006). More generally we also have:

**Theorem.2.6** For a *d.f* of the form $H(x) = \dfrac{1}{1+\psi(x)}$, $\dfrac{1}{1+a\psi(x)}$ is a *d.f* of a geometric($p$)-max for any $a>0$.

*Proof.* $\dfrac{1}{1+a\psi(x)} = \dfrac{1}{a}\left(\dfrac{1}{a}+\psi(x)\right)^{-1} = \dfrac{1}{a}\left[\dfrac{1}{a}+\left(\dfrac{1}{H(x)}-1\right)\right]^{-1}$

$$= \dfrac{1}{a}H(x)\left[\dfrac{H(x)+a-aH(x)}{a}\right]^{-1}$$

$$= \dfrac{1}{a}H(x)\left[1-\left(1-\dfrac{1}{a}\right)H(x)\right]^{-1}$$

$$= \sum_{k=1}^{\infty} \frac{1}{a} \left(1 - \frac{1}{a}\right)^{k-1} [H(x)]^k$$

$$= P\big(Max(X_1, X_2, \ldots X_{N(p)}) \leq x\big), \ p = \tfrac{1}{a},$$

where $N(p)$ is a geometric($p$) $r.v$ and $X_i$'s are $i.i.d$ with $d.f$ $H(x)$ proving the assertion.

This is the max-analogue of lemma.3.1 in Pillai (1990). In particular we have:

**Theorem.2.7** Geometric($p$)-max of $i.i.d$ G-gamma-MID $(\beta)$ variables is G-gamma-MID $\left(\frac{\beta}{p}\right)$ for any $p \in (0,1)$ where the geometric($p$) $r.v$ is independent of the components.

*Proof.* The $d.f$ of geometric($p$)-max of $i.i.d$ G-gamma-MID $(\beta)$ variables is given by;

$$\frac{p/\{1 + \beta \log(1 + \psi(x))}{1 - (1-p)/\{1 + \beta \log(1 + \psi(x))\}}$$

$$= \frac{p}{p + \beta \log(1 + \psi(x))} = \frac{1}{1 + \frac{\beta}{p} \log(1 + \psi(x))},$$

which is the $d.f$ of a G-gamma-MID $\left(\frac{\beta}{p}\right)$ law.

**Remark.2.2** The invariance property under geometric($p$)-max described above is new and is different from the invariance property under geometric($p$)-max up-to a scale-change in theorem.2.5 characterizing geometric-max-semi-stable laws. In theorem.2.7 it is invariance up-to a change of shape parameter.

**3. Processes related to gamma-MID and G-gamma-MID laws.**

**3.1 Random Time Changed Extremal Processes.** Extremal processes (EP) are processes with increasing right continuous sample paths and independent max-increments. The univariate marginals of an EP determine its finite dimensional distributions. Also, the max-increments of an EP are MID. The EPs will be referred to by the distribution of their max-increments. Pancheva, *et al.* (2006) has discussed random time changed or compound EPs and their theorem.3.1 together with property.3.2 reads: Let $\{Y(t), t \geq 0\}$ be an EP having homogeneous max-increments with $d.f$ $F_1(y) = exp\{-t\mu([\lambda, y]^c)\}$, $y \geq \lambda$, $\lambda$ being the bottom of the rectangle $\{F > 0\}$ and $\mu$ the exponential measure of $Y(1)$, that is, $\mu([\lambda, y]^c) = -\log F(y)$. Let $\{T(t), t \geq 0\}$ be a non-negative process independent of $Y(t)$ having stationary, independent and additive increments with Laplace transform (LT) $\phi$. If $\{X(t), t \geq 0\}$ is the compound EP obtained by randomizing the time parameter of $Y(t)$ by $T(t)$, then $X(t) = Y(T(t))$ and its $d.f$ is:

$$P\{X(t)<x\} = \{\varphi(\mu([\lambda,x]^c))\}^t.$$

Pancheva, *et al.* (2006) also showed that in the above setup $Y(T(t))$ is also an EP. Extending the discussion in Pancheva, *et al.* (2006) Satheesh and Sandhya (2006) showed that the EP obtained from a random time changed (by compounding a) max-semi-stable EP is max-semi-selfdecomposable(*b*) if the compounding process is semi-selfdecomposable. Here we have some more results in this direction.

**Theorem.3.1** The EP obtained by compounding a gamma-MID$(\beta)$ EP having homogeneous max-increments is G-gamma-MID$(\beta)$ if the compounding process is unit exponential.

*Proof.* If $\{Y(t)\}$ is a gamma-MID$(\beta)$ EP and $\{T(t)\}$ is unit exponential with *d.f G* then:

$$\frac{1}{1+\beta \log(1+\psi(x))} = \int_0^\infty e^{-t\beta \log(1+\psi(x))} dG(t), \text{ which proves the assertion.}$$

**Theorem.3.2** The EP $\{Y(T(t))\}$ obtained from a random time changed EP $\{Y(t)\}$ having homogeneous max-increments with *d.f* $e^{-H(x)}$ is G-gamma-MID$(\beta)$ if the compounding process $\{T(t)\}$ is G-gamma$(\beta)$.

*Proof.* We know that the LT of a G-gamma$(\beta)$ law is $\dfrac{1}{1+\log(1+\lambda)^\beta}$. If $G$ denote its *d.f* then;

$$\frac{1}{1+\beta \log(1+\psi(x))} = \int_0^\infty e^{-t\psi(x)} dG(t), \text{ proving the assertion.}$$

**3.2. An Auto-regressive Model.** Now consider a first order max-autoregressive (max-AR(1)) model described as below. Here a sequence of *r.v*s $\{X_n, n>0 \text{ integer}\}$ defines the max-AR(1) scheme if for some $0<p<1$ there exists an innovation sequence $\{\varepsilon_n\}$ of *i.i.d r.v*s such that;

$$
\left.
\begin{array}{l}
X_n = X_{n-1}, \text{ with probability } p \\[2mm]
\phantom{X_n} = X_{n-1} \vee \varepsilon_n, \text{ with probability } (1\text{-}p).
\end{array}
\right\} \tag{1}
$$

In terms of *d.f*s and assuming stationarity this is equivalent to;

$F(x) = F(x)\{p + (1\text{-}p) F_\varepsilon(x)\}$. That is;

$$F(x) = \frac{pF_\varepsilon(x)}{1-(1-p)F_\varepsilon(x)}$$

Hence $\{X_n\}$ is a geometric($p$)-max of the innovation sequence $\{\varepsilon_n\}$. Invoking theorem.2.7 we have proved;

**Theorem.3.3** In the max-AR(1) structure (1) the sequence $\{X_n\}$ and the innovation sequence $\{\varepsilon_n\}$ are related as follows for any $p \in (0,1)$. $\{X_n\}$ is G-gamma-MID$(\beta)$ variables *iff* $\{\varepsilon_n\}$ is G-gamma-MID$(\frac{\beta}{p})$. Equivalently,

**Theorem.3.4** A necessary and sufficient condition for an max-AR(1) process $\{X_n\}$ with the structure in (1) is stationary Markovian for any $p \in (0,1)$ with G-gamma-MID$(\beta)$ distribution is that the innovation's are G-gamma-MID$(\frac{\beta}{p})$ distributed.

**References.**

Balkema, A. A and Resnick, S. (1977). Max-infinite divisibility, *J. Appl. Probab*. **14**, 309-319.

Kaufman, E (2001). *Statistical Analysis of Extreme Values, From Insurance, Finance, Hydrology and Other Fields*, Extended 2[nd] Edition, Birkhauser.

Mittnick, S. and Rachev, S. T (1993). Modelling asset returns with alternative stable distributions, *Econ. Rev.*, **12**, 161-330.

Mohan, N. R (1998). On geometrically max-infinitely divisible laws, *J. Ind. Statist. Assoc.* **36**, 1-12.

Pancheva, E; Kolkovska, E. T and Jordanova, P. K (2006). Random time changed external processes, *Probab. Theory Appl*. **51**, 752-772.

Pillai, R. N (1990). Harmonic mixtures and geometric infinite divisibility, *J. Ind. Statist. Assoc.* **28**, 87-98.

Rachev, S. T (1993). Rate of convergence for maxima of random arrays with applications to stock returns, *Statist. and Decisions*, **11**, 279-288.

Rachev, S. T and Resnick, S (1991). Max-geometric infinite divisibility and stability, *Comm. Statist. – Stoch. Mod.* **7**, 191-218.

Sandhya, E and Satheesh, S (2007). Geometric generalized gamma distributions and related processes, *to appear in Appl. Math. Sci.*, **2008**.

Satheesh, S (2002). Aspects of randomization in infinitely divisible and max-infinitely divisible laws, *ProbStat Models* **1**, 7-16.

Satheesh, S and Sandhya, E (2006*b*). Max-semi-selfdecomposable laws and related processes, *Statist. Probab. Lett.* **76**, 1435-1440.